\documentclass[10pt]{article}
\usepackage{amsmath,amssymb,pictex}
\usepackage[margin=20pt,font=small,labelfont=bf]{caption}

\newcommand{\nd}{\mathop{\hbox{\rm nd}}}

\begin{document}
\thispagestyle{empty}

\title{\bfseries A New Lower Bound for the Domination Number of
  Complete Cylindrical Grid Graphs}
\author{%
    David R. Guichard\\%
    Whitman College\\%
}
\date{}
\maketitle

\begin{abstract} 
We use a dynamic programming algorithm to
establish a lower bound on the domination number 
 of complete grid graphs of the form $C_n\square P_m$, that is, the
 Cartesian product of a cycle $C_n$ and a
 path $P_m$, for $m$ and $n$ sufficiently large.

\vspace{.1in}
\noindent keywords: cylindrical grid graph, domination number

\vspace{.1in}
\noindent AMS classification: 05C69
\end{abstract}

\section{Introduction} A set $S$ of vertices in a graph $G=(V,E)$ is
called a {\it dominating set\/} if every vertex $v\in V$ is either in
$S$ or adjacent to a vertex in $S$.  The domination number of $G$,
$\gamma(G)$, is the minimum size of a dominating set.

Let $P_m$ denote the path on $m$ vertices and $C_n$ the cycle on $n$
vertices; the {\it complete cylindrical grid graph\/} or
{\it cylinder} is the product
$C_n\square P_m$. That is, if we denote the vertices of $C_n$ by
$u_1,u_2,\ldots,u_n$ and the vertices of $P_m$ by $w_1,\ldots,w_m$,
then $C_n\square P_m$ is the graph with vertices $v_{i,j}$, $1\le
i\le n$, $1\le j\le m$, and $v_{i,j}$ adjacent to $v_{k,l}$ if
$i=k$ and $w_j$ is adjacent to $w_l$ or if
$j=l$ and $u_i$ is adjacent to $u_k$.
It will  be useful to think of this graph as $P_n\square P_m$,
with the edge paths of length $m$ glued together, that is, connected
with new edges.

P. Pavli{\v c} and J. {\v Z}erovnik\cite{pavlic-zerovnik:dom-no-cyl-graphs}
established upper bounds for
the domination number of $C_n\square P_m$, and Jos\'e Juan Carre{\~n}o
et al.\cite{carreno-et-al:lower-bound-dom-no-cyl-graphs} established
non-trivial lower bounds. For $n\equiv 0\pmod 5$ the bounds agree, so
the domination number is known exactly. Here we improve the lower
bounds, except of course in the case that $n\equiv 0\pmod 5$. The
method is similar, based on a technique first used in
Guichard\cite{drg:dom-grid-digraph-lower-bound}, and later in
Gon\c{c}alves, et al.\cite{goncalves-et-al:domination_number_of_grids},
but we use a different programming technique than that of
\cite{carreno-et-al:lower-bound-dom-no-cyl-graphs}.

\section{Getting a lower bound} 
\begin{figure}[htb]
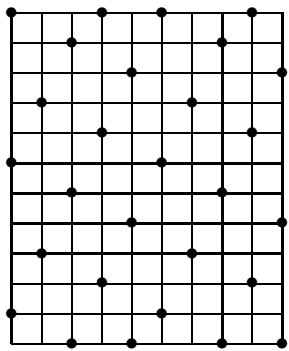

\hbox to \hsize{\hss\beginpicture
\setcoordinatesystem units <4mm,4mm>
\setplotarea x from 0 to 9, y from 0 to 11
\grid 9 11
\multiput {$\bullet$} at 2 0 4 0 7 0 9 0 /
\multiput {$\bullet$} at 0 1 5 1 /
\multiput {$\bullet$} at 3 2 8 2 /
\multiput {$\bullet$} at 1 3 6 3 /
\multiput {$\bullet$} at 4 4 9 4 /
\multiput {$\bullet$} at 2 5 7 5 /
\multiput {$\bullet$} at 0 6 5 6 /
\multiput {$\bullet$} at 3 7 8 7 /
\multiput {$\bullet$} at 1 8 6 8 /
\multiput {$\bullet$} at 4 9 9 9 /
\multiput {$\bullet$} at 2 10 7 10 /
\multiput {$\bullet$} at 0 11 3 11 5 11 8 11 /
\endpicture\hss}
\caption{The cylinder $C_{10}\square P_{12}$ has
domination number 28. (Vertices on the right side are adjacent to the
corresponding vertices on the left side.)}
\label{fig:12x10}
\end{figure}

A vertex in $C_n\square P_m$ dominates at most five vertices,
including itself, so certainly $\gamma(C_n\square P_m)\ge nm/5$. If we
could keep the sets dominated by individual vertices from overlapping,
we could get a dominating set with approximately $nm/5$ vertices, and
indeed we can arrange this for much of the graph, with the exception
of the the top and bottom copies of $C_n$ in which the vertices have
only 3 neighbors, and, except when $n\equiv 0\pmod 5$, in the leftmost
and rightmost columns of $P_n\square P_m$ where each vertex in the
leftmost column is adjacent to the corresponding vertex in the
rightmost column.  Figure~\ref{fig:12x10} shows one of the nice
examples, when $n$ is divisible by 5.

Suppose $S$ is a subset of the vertices of $C_n\square P_m$. Let $N[S]$ be
the set of vertices that are either in $S$ or adjacent to a member of
$S$, that is, the vertices dominated by $S$.  Define the {\it wasted
domination\/} of $S$ as $w(S)=5|S|-|N[S]|$, that is, the number of
vertices we could dominate with $|S|$ vertices in the best case, less
the number actually dominated. When $S$ is a dominating set,
$|N[S]|=mn$, and if $w(S)\ge L$ then $|S|\ge (L+mn)/5$. Our goal now
is to find a lower bound $L$ for $w(S)$.

\begin{figure}[tb]
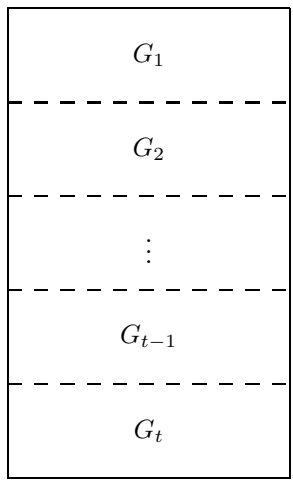

\hbox to \hsize{\hss
\beginpicture
\setcoordinatesystem units <2.5mm,2.5mm>
\setplotarea x from 0 to 15, y from 0 to 25
\axis left /
\axis right /
\axis top /
\axis bottom /
\setdashes
\putrule from 0 5 to 15 5
\putrule from 0 10 to 15 10
\putrule from 0 15 to 15 15
\putrule from 0 20 to 15 20
\putrule from 0 25 to 15 25
\put {$G_1$} at 7.5 22.5
\put {$G_2$} at 7.5 17.5
\put {$\vdots$} at 7.5 12.5
\put {$G_{t-1}$} at 7.5 7.5
\put {$G_t$} at 7.5 2.5
\endpicture\hss}
\caption{Partitioned cylinder.}
\label{fig:partitioned}
\end{figure}

Suppose a cylinder $C_n\square P_m$ is partitioned into
subgraphs as indicated in Figure~\ref{fig:partitioned}, where each $G_i$
is a subgraph $C_n\square P_{m_i}$. 
Let $S$ be a
dominating set for $G$ and $S_i=S\cap V(G_i)$.  Then
\begin{equation}\label{eq:fundamental inequality}
w(S)\ge \sum_{i=1}^t w(S_i).
\end{equation}
Note that in computing $w(S_i)$ we consider $S_i$ to be a subset of $V(G)$,
not of $V(G_i)$ (this affects the computation of $N[S_i]$). To verify
the inequality, note that the following inequalities are equivalent:
\begin{align*}
  w(S)&\ge \sum_{i=1}^t w(S_i)\\
  5|S|-|N[S]|&\ge \sum_{i=1}^t (5|S_i|-|N[S_i]|)\\
  5|S|-|N[S]|&\ge \sum_{i=1}^t 5|S_i|-\sum_{i=1}^t |N[S_i]|\\
  |N[S]| &\le \sum_{i=1}^t |N[S_i]|.  
\end{align*}
The last inequality is satisfied, since each vertex in $N[S]$ is
counted at least once by the expression on the right. 

Note that $S_i$ is a set that dominates all the vertices of $G_i$
except possibly some vertices in the top or bottom row of $G_i$ (or in
the cases of $G_1$ and $G_t$, in the bottom row and top row,
respectively).  Let us say that a set that dominates a
cylinder $G$, with the exception of some vertices on the
top or bottom edges, {\it almost dominates\/} $G$.  Given a 
cylinder $H=C_n\square P_{m_i}$ (namely, one of the $G_i$),
What we want to know is the value of
\begin{equation}\label{eq:desired minimum}
  \min_A w(A),
\end{equation}
taking the minimum over sets $A$ that almost dominate
$H$ and computing $w(A)$ as if $A$ were a subset of a larger graph
$C_n\square P_{m_i+2}$ in which $H$ occupies the middle $m_i$ rows, or in
the case of $G_1$ or $G_t$, $A$ is a subset of $C_n\square P_{m_i+1}$ in
which $H$ occupies the top $m_i$ rows.  If we can compute this minimum
for (small) fixed $m_i$ and any $n$, we can choose $G_1$ through $G_t$
with a small number of rows and get lower bounds on $w(S_i)$ for any
dominating set $S$ of the original $C_n\square P_m$.

\section{The algorithm}

We describe the algorithm for $G_1$ and $G_t$ (which of course are
isomorphic); the algorithm for the other graphs $G_i$ is nearly
identical, and we describe it more briefly.  Imagine a
cylinder $C_n\square P_m$ with a designated subset $S$
of the vertices. Recall that the vertices are denoted by $v_{i,j}$,
$1\le i\le n$, $1\le j\le m$ (say, numbering left to right and bottom
to top).  We describe a column, say column number $i$, in such a
diagram by a state vector $\bf s$, in which ${\bf s}_j$ is $0$ if
vertex $v_{i,j}$ is in $S$, $1$ if vertex $v_{i,j}$ is adjacent to a
member of $S$ in column $i$ or column $i-1$, and $2$ otherwise. For
example, the second column from the right in Figure~\ref{fig:12x10}
has state vector $(1,1,0,1,2,1,1,0,1,2,1,0)$.  Let $|{\bf s}|$ denote
the number of zeros in $\bf s$.

Given a state vector $\bf s$, we append a column 0 at the left of
$P_n\square P_m$. Let $X$ be the set of vertices in this column
corresponding to the 0 entries in $\bf s$, and let $Y$ be the set of
vertices corresponding to the 2 entries in $\bf s$.  An
{\it $({\bf s,t})$-almost-domination\/} of $P_n\square P_m$ is a subset $S$ of
the vertices such that $X\cup S$ dominates the first $n-1$ columns of
$P_n\square P_m$ and the elements of $Y$,
except possibly vertices in the first (i.e., bottom)
row, and for which the state vector of the final column is $\bf t$.

Suppose $S$ is a subset of the vertices of $P_i\square P_j$ and denote
by $w_{i,j}(S)$ the value of $w(S)$ computed in $P_{i+1}\square
P_{j+1}$, in which $P_i\square P_j$ occupies the top $j$ rows and
leftmost $i$ columns. Let
$$w_{i,j}({\bf s,t})=\min_S w_{i,j}(S),$$
taking the minimum over all $({\bf s,t})$-almost-dominations of
$P_i\square P_j$.
If there is no $({\bf s,t})$-almost-domination of
$G_{i,j}$, let $w_{i,j}({\bf s})=\infty$.

Finally, to compute the desired minimum (equation~\ref{eq:desired
  minimum}), we compute
$$\min_{\bf s} w_{i,j}({\bf s,s}),$$
since an $({\bf s,s})$-almost-domination of $P_i\square P_j$ almost
dominates $C_i\square P_j$.

Let ${\cal P}({\bf t})$ be the set of state vectors $\bf u$ such that $\bf
u$ is the state vector of the next to last column in an 
$({\bf s,t})$-almost-domination of $P_n\square P_m$.
Then
$$w_{n,m}({\bf s,t})=\min_{{\bf u}\in{\cal P}({\bf t})} 
   \left( 5|{\bf t}| - \nd({\bf u},{\bf t})+w_{n-1,m}({\bf s,u}) \right),
$$
where $\nd({\bf u},{\bf t})$, the number of newly dominated vertices,
may be computed as follows.

\newcounter{ctr}

\begin{list}{\arabic{ctr}.}{\usecounter{ctr}}

\item $\nd = 0$
\item For each $j=1,\dots,m$ for which ${\bf t}_j=0$ and
${\bf u}_j=2$, add 1 to $\nd$. This counts the newly dominated vertices
$v_{n-1,j}$.
\item For each $j=1,\dots,m$ for which ${\bf t}_j\le 1$ and
${\bf u}_j\ge 1$, add 1 to $\nd$. This counts the newly
dominated vertices $v_{n,j}$.
\item For each $j=1,\dots,m$ for which ${\bf t}_j=0$, add 1 to
$\nd$. This counts the newly dominated vertices $v_{n+1,j}$.
\item If ${\bf t}_1=0$, add 1 to $\nd$. This counts the newly
  dominated vertex below vertex $v_{n,1}$, recalling that we compute
  $w(S)$ in $P_n\square P_m$ with an extra bottom row.
\end{list}

\noindent
Now, given some $n$, the algorithm to compute $w_{n,m}({\bf s},{\bf
  t})$, $i=1,\dots,n$, is:

\begin{list}{\arabic{ctr}.}{\usecounter{ctr}}

\item {\bf Initialization.} Set $w_{0,m}({\bf s},{\bf u})=0$ if
${\bf u}= {\bf s}$, and $\infty$ otherwise.

\item {\bf Iteration.} Suppose that $i\le n$ and that
$w_{i-1,m}({\bf s},{\bf u})$ has been computed for all ${\bf u}$.
Then for each $\bf t$, set 
$$ w_{i,m}({\bf s},{\bf t})=\min_{{\bf u}\in{\cal P}({\bf t})} 
         \left( 5|{\bf t}| - \nd({\bf u},{\bf t})+
           w_{i-1,m}({\bf s},{\bf u}) \right).
$$ 

\end{list}

Thus, for fixed $m$ and any $n$, we can compute $\min_{\bf s}
w_{n,m}({\bf s},{\bf s})$, by computing $w_{i,m}({\bf s},{\bf t})$ for
all ${\bf s}$, ${\bf t}$, and $1\le i\le n$.
Of course, what we want is to know this
value for any $n$ without an infinite amount of work.  Livingston and
Stout~\cite{ls:constant-time-dominating-sets} and
Fisher~\cite{dcf:domination-grid-graphs} independently thought of
looking for a sort of periodicity in the values of $\gamma(P_n\square
P_m)$ for fixed $m$. Since they succeeded, we might hope that for
fixed $m$, there are $N$, $p$, and $q$ so that for $n\ge N$ and all
$\bf s$ and $\bf t$,
$$w_{n,m}({\bf s},{\bf t})=w_{n-p,m}({\bf s},{\bf t} )+q.$$
In this case, after a finite amount of computation, we could determine 
$\min_{\bf s} w_{n,m}({\bf s},{\bf s})$ for all $n$.

It is easy to modify the algorithm so to check for this
periodicity. When we do this, we find that for $n\ge 65$, 
$$\min_{\bf s} w_{n,10}({\bf s},{\bf s}) =
\min_{\bf s} w_{n-1,10}({\bf s},{\bf s})+1 =
n.$$
Thus, for $m\ge 20$ and $n\ge 64$,  if $S$ is a dominating set
in $C_n\square P_m$,
$$w(S)\ge \sum_{k=1}^t w(S_k) \ge w(S_1)+w(S_t) \ge 2n,$$
using the inequality (\ref{eq:fundamental inequality}), and so
$$|S|\ge (mn + 2n)/5.$$
Jos\'e Juan Carre{\~n}o 
et al.\cite{carreno-et-al:lower-bound-dom-no-cyl-graphs} have independently
arrived at the same conclusion, using a substantially different
algorithm. When $n\cong 0\pmod{5}$, $(mn + 2n)/5$ is also known to be
an upper bound, so that $\gamma(C_n\square P_m)=(mn + 2n)/5$ (in fact,
this is known to be correct for $n\ge 5$).
The
implication, of course, is that for optimal $S$, 
$w(S_k)=0$ for $1<k<t$, when $n\cong 0\pmod{5}$.
This is not true in general, so we improve our lower bound by
computing a lower bound on $w(S_k)$, $1<k<t$.

The only change required is to redefine an $({\bf s,t})$-almost-domination as follows:
Given a state vector $\bf s$, we append a column 0 at the left of
$P_n\square P_m$. Let $X$ be the set of vertices in this column
corresponding to the 0 entries in $\bf s$, and let $Y$ be the set of
vertices corresponding to the 2 entries in $\bf s$.  An
$({\bf s,t})$-almost-domination of $P_n\square P_m$ is a subset $S$ of
the vertices such that $X\cup S$ dominates the first $n-1$ columns of
$P_n\square P_m$ and the elements of $Y$,
except possibly vertices in the top and bottom rows,
and for which the state vector of the final column is $\bf t$.
Corresponding to this change, in the computation of $nd$, we add a
sixth step:
\begin{list}{\arabic{ctr}.}{\usecounter{ctr}}{\setcounter{ctr}{5}}
\item If ${\bf t}_m=0$, add 1 to $\nd$. This counts the newly dominated
vertex above $v_{n,m}$, recalling that we compute $w(S)$ in
$C_n\square P_m$ as if it occupies the middle $m$ rows of a copy of
$C_n\square P_{m+2}$.
\end{list}

Proceeding as before, we find that $\min_{\bf s} w_{n,10} ({\bf s},{\bf s})
=\min_{\bf s} w_{n-5,10} ({\bf s},{\bf s})$, when $n\ge
12$. Specifically, we find that
$\min_{\bf s} w_{n,10}({\bf s},{\bf s})$ is 0, 6, 5, 9,
or 6 as $n$ is 0, 1, 2, 3, or 4 $\pmod 5$.
Thus, with $a$ equal to 0, 6, 5, 9, or 6 as appropriate, we find that
$$|S|\ge {1\over 5}((m+2)n+\lfloor {m-20\over 10}\rfloor\cdot a).$$
That is, lower bounds for the domination number of $C_n\square P_m$,
when $m\ge 20$ and $n\ge 64$, 
are:
\begin{align*}
  {(m+2)n\over 5},&\qquad n\equiv 0\pmod 5\\
  {(m+2)n\over 5}+{6\over 5}\lfloor {m-20\over 10}\rfloor,&\qquad n\equiv 1\pmod 5\\
  {(m+2)n\over 5}+\lfloor {m-20\over 10}\rfloor,&\qquad n\equiv 2\pmod 5\\
  {(m+2)n\over 5}+{9\over 5}\lfloor {m-20\over 10}\rfloor,&\qquad n\equiv 3\pmod 5\\
  {(m+2)n\over 5}+{6\over 5}\lfloor {m-20\over 10}\rfloor,&\qquad n\equiv 4\pmod 5.
\end{align*}
Known upper bounds (see \cite{pavlic-zerovnik:dom-no-cyl-graphs})
for the domination number of $C_n\square P_m$ are:
\begin{align*}
  {(m+2)n\over 5},&\qquad n\equiv 0\pmod 5\\
  {(m+2)n\over 5}+{7\over 40}(m+2),&\qquad n\equiv 1\pmod 5\\
  {(m+2)n\over 5}+{1\over 10}(m+2),&\qquad n\equiv 2\pmod 5\\
  {(m+2)n\over 5}+{2\over 5}(m+2),&\qquad n\equiv 3\pmod 5\\
  {(m+2)n\over 5}+{1\over 5}(m+2),&\qquad n\equiv 4\pmod 5.
\end{align*}

For $n\equiv 2\pmod 5$ the lower and upper bounds are quite close, but
for the other non-zero values of $n\bmod 5$ there is considerable room
for improvement. It seems likely that the upper bounds are closer to
the true values, as our computation allows vertices on the boundary
(that is, the top and bottom rows) of
the subgraphs $G_k$ to remain undominated. A small increase in the
value of $a$ in each case would eliminate most of the gap.

When $m\bmod 10$ is non-zero, we have effectively ignored one of the
$G_i$, that is, used zero as a lower bound for one of the $w(S_i)$.
We can improve our lower bounds very slightly (by a small constant) by
correcting this. For example, for $m>20$ and 
$m\equiv 8\pmod{10}$, we can let all
but one of the $G_i$ have height 10, and the remaining (interior)
graph, say $G_2$, have height 8. Then we run the algorithm again for
height 8 graphs. While we have in fact done the additional
computations, the improvement is very slight, so we omit the results.

Our approach gives us lower bounds for $m\ge 20$;
Crevals~\cite{crevals:domination_cylinder_graphs} computes exact
values for $m\le 22$ and all $n$. He also computes exact values for
$n\le 30$ and all $m$. In the course of our computations, we also
obtain lower bounds for $12\le n<64$ (with only $n>30$ of interest due
to the Crevals results), but they do not seem sufficiently
illuminating to include here.

\bibliography{guichard}
\bibliographystyle{abbrv}

\end{document}